\documentclass[graybox]{svmult}

\usepackage{graphicx} %
\usepackage{amsmath}
\usepackage{amstext}
\usepackage{amsfonts}
\usepackage{amssymb}
\usepackage{color}

\usepackage{longtable}
\usepackage{bm}
\usepackage{float}

\newcommand{\be}{\begin{equation}}
\newcommand{\ee}{\end{equation}}
\newcommand{\bea}{\begin{eqnarray}}
\newcommand{\eea}{\end{eqnarray}}
\newcommand{\Eq}[1]{Eq.\,(\ref{#1})}

\newcommand{\eqa}{\begin{equation}}
\newcommand{\eqz}{\end{equation}}
\newcommand{\eqma}{\begin{eqnarray}}
\newcommand{\eqmz}{\end{eqnarray}}

\newcommand{\pot}[1]{\ensuremath{^{\left(#1\right)}}}

\newcommand{\kl}[1]{\ensuremath{\left(#1\right)}}
\newcommand{\kll}[1]{\ensuremath{\left\{#1\right\}}}
\newcommand{\klll}[1]{\ensuremath{\left[#1\right]}}

\newcommand{\norm}[1]{\left\lVert#1\right\rVert}
\usepackage{array}
\newcolumntype{R}[1]{>{\raggedleft  \arraybackslash}p{#1}@{} }
\newcolumntype{C}[1]{>{\centering \arraybackslash}p{#1}@{} }

\usepackage{mathptmx}        %
\usepackage{helvet}          %
\usepackage{courier}         %
\usepackage{type1cm}         %

 \usepackage{makeidx}         %
\usepackage{graphicx}        %
\usepackage{multicol}        %
\usepackage[bottom]{footmisc}%

\usepackage{amsmath,amssymb,bbm}
\usepackage{graphicx}
\usepackage[utf8]{inputenc}
\usepackage{url}
\usepackage{microtype}
\usepackage{import}
\usepackage{color}
\usepackage{epic}
\usepackage{tikz}
\usepackage{etex}
\usepackage{subfig}

\begin{document}

\title*{Efficiently transforming from  values of a function on a sparse grid to basis coefficients}
\titlerunning{Basis coefficients from function values  }   %
\author{Robert Wodraszka \and Tucker Carrington Jr.}

\institute{Chemistry Department, Queen's University, Kingston, Ontario K7L 3N6, Canada}

\maketitle

\abstract

In many contexts it is necessary to determine coefficients of a basis expansion of a function $
{f}\kl{x_1, \ldots, x_D} $ from values of the function at points on a sparse grid.   Knowing the coefficients, one has an interpolant or a surrogate.  For example, such 
coefficients are used in uncertainty quantification.  In this chapter, we present an efficient method for computing the coefficients.   
It uses basis functions that, like the familiar piecewise linear hierarchical functions, are zero at points in previous levels.  
They  are linear combinations of \textit{any}, e.g. global,  nested basis functions $\varphi_{i_k}\pot{k}\kl{x_k}$.
Most importantly, the transformation from function values to basis coefficients is done, exploiting the nesting, by evaluating sums sequentially.  
When the number of functions in level $ \ell_k $ equals $ \ell_k $  %
(i.e. when the level index is increased by one, only one point (function) is added)  %
and the basis function  indices satisfy  $ {\norm{\mathbf{i}-\mathbf{1}}_1 \le b}$,
 the cost of the transformation scales as  
 $\mathcal{O}\kl{D \klll{\frac{b}{D+1} + 1} N_\mathrm{sparse}}$,
where $  N_\mathrm{sparse} $  is the number of points on the sparse grid.  %
 We compare the cost of doing the transformation with sequential sums to the cost of other methods in the literature.

\section{Introduction}
Sparse grids are often used to interpolate multi-dimensional functions. \cite{BaNoRi,KlWo,PfBook}
An interpolant for a function $f\kl{x_1, x_2, \ldots, x_D}$ is usually built from 
a set of basis functions and a set of interpolation points. In this chapter,
the basis functions are products of 1-D functions which are denoted 
$\varphi_{i_k}\pot{k}\kl{x_k}$, where $i_k = 1, 2, \ldots, n_k$ and $k = 1, 2, \ldots, D$.  
$ D $ is the number of dimensions and $ n_k $ is the number of basis functions for dimension $ k. $
All interpolants in this chapter can be written in the general form  
\begin{align}
   \bar{f}\kl{x_1, \ldots, x_D} &= \sum_{{\mathbf{i}} } C_\mathbf{i} \;
   \prod_{k=1}^D \varphi_{i_k}\pot{k}\kl{x_k},
   \label{inter}
\end{align}
where the sum is over a set of  $\mathbf{i}=\kl{i_1, \ldots, i_D}$ indices.   
Throughout this chapter:  a superscript in brackets indicates a particular dimension;  
$\varphi_{j_k}\pot{k}$ is the $j_k$-th basis function for dimension  $ k $;   
vectors, matrices, and tensors are in bold type, level indices are  script $ \ell $; 
the lower case letter $ a $ labels a point; 
and basis set functions are labelled by $ i,j,i',j' $.  
The coefficients $C_\mathbf{i}$ are determined so that 
the interpolant is equal to known function values at the interpolation points, 
$\kl{ r_{a_1}\pot{1}, r_{a_2}\pot{2},  \ldots, r_{a_D}\pot{D}}        $.  %
Note that $ r_{a_k}\pot{k} $ is a point in dimension $k$ and that $ a_k  $  labels a point. %
We shall assume throughout the chapter that the number of 
interpolation points is equal to the number of basis functions.

A simple and general interpolation method is obtained by using special 1-D basis functions, $G_{a_k}\pot{k}\kl{x_k}$,
called Lagrange type functions, that are equal to zero at all %
points except one,
a tensor product basis, and a tensor product grid of interpolation points.
In this case, the coefficients $C_\mathbf{i}$  are simply the known values of the function at the points.   
The most common Lagrange type functions are the Lagrange polynomials  that span the same polynomial space as the monomial basis $\kll{x_k^0, x_k^1,  \ldots, x_k^{n-1}}$.   
However, it is possible to make Lagrange type functions that span the space spanned by any 1-D basis, see \cite{AvCaV,MaOcc}, %
\begin{align}  
   G_{a_k}\pot{k}\kl{x_k} &= \sum_{j_k=1}^{n_k}
   [(\mathbf{B\pot{k}   })^{-1}]_{j_k,a_k}
    \varphi_{j_k}\pot{k}\kl{x_k},~a_k = 1, 2, \ldots, n_k ~.   
   \label{LTF}  
\end{align}
$ n_k $ is the number of basis functions for dimension $ k $.
A Lagrange type function   is   labelled by the point at which it is centred.  In this chapter, points for dimension $ k $ are labelled
 by  $ a_k $ or $ a'_k $.  %
 In \Eq{LTF},   $(\mathbf{B\pot{k}})_{a_k,j_k} = 
  \varphi_{j_k}\pot{k}\kl{    r_{a_k}\pot{k} }$.    %
Making these Lagrange type functions requires inverting a small $n_k \times n_k$ matrix.  
Instead of Lagrange type functions, it is also common to use a basis composed of piecewise linear functions
(hat functions), called a "nodal" basis. %
\cite{BuGr} %
 In this nodal basis, the 
coefficients $C_\mathbf{i}$ are also the known values of the function at the points.

Any tensor product basis is afflicted by the curse of dimensionality.  
A better multi-dimensional basis can be made by using sparse grid or Smolyak type ideas and 
1-D basis functions that are importance-ordered. 
\cite{BuGr}   %
An importance-ordered 1-D basis,  $\varphi_{i_k}\pot{k}\kl{x_k}$, 
is one in which a basis function is more important if its  value of $i_k$ is smaller.
An important basis function is one whose coefficient is large. 
This is the famous idea of Archimedes. \cite{Archi}  %
For the purpose of interpolating with the sparse-grid Ansatz,  it is best to use ZAPPL functions (see section 2), made from an importance-ordered basis.
To interpolate with Smolyak's idea one frequently uses \cite{BaNoRi}     
\begin{align}  %
I\kl{D,b}
   &= \sum_{\norm{\bm{\ell}-\mathbf{1}}_1 \le b}    %
   \Delta^{\ell_1} \otimes \Delta^{\ell_2} \ldots  \otimes   \Delta^{\ell_D}~,
   \label{bart}
\end{align}
where $\ell_k$ labels a level. 
\begin{align}
   \Delta^{\ell_k}  &=     U^{\ell_k} -  U^{\ell_k-1} 
\end{align}
and $U^{\ell_k}$ is a 1-D interpolation rule;  $U^0 = 0$.  \cite{BaNoRi}  
In the restriction on the sum,    %
$\norm{\bm{\ell}-\mathbf{1}}_1=\sum_{k=1}^D \left|l_k - 1\right|$,  where %
$l_k=1, \ldots, b + 1 \; \forall k=1, \ldots, D$.  %
The importance-ordered basis is divided into levels and for each level there is a corresponding set of points.   
In level $\ell_k$, for coordinate $k$, there are $m_k\kl{\ell_k}$ points and basis functions.   
In this chapter, for simplicity, we shall set $m_k\kl{\ell_k} = \ell_k$, but the same ideas can be implemented 
when $ m_k\kl{\ell_k}  \ge \ell_k$. \cite{AvCaV} %
Everywhere in this chapter, we shall assume that the sequences of points are nested, 
i.e., the set of points with $m_k\kl{\ell_k}$ points includes all the points in the set with $m_k\kl{\ell_k-1}$  %
In this chapter, the grids included in the sum in \Eq{bart} are those that satisfy the condition 
$\norm{\bm{\ell}-\mathbf{1}}_1 \le b  $, but other choices are possible.   \cite{AvCaV}  %
The space spanned by the pruned (restricted)
basis is smaller than the space spanned by the  tensor product basis.  
If both the function being interpolated and the basis functions are smooth, the pruning is effective.   
In general, formulations that obviate the sum over levels are less costly than   \Eq{bart}   which requires a sum over levels.

The interpolant made from  \Eq{bart} can be written in terms of Lagrange type functions or in terms of the functions
$\varphi_{i_k}\pot{k}\kl{x_k}$, from which the Lagrange type functions are made.  %
If written in terms of $\varphi_{i_k}\pot{k}\kl{x_k}$, the interpolant is 
\begin{eqnarray} 
\klll{I(D,b) ~ f}\kl{x_1, \ldots, x_D} =  %
   \sum_{\norm{\mathbf{i}-\mathbf{1}}_1 \le b}  
   C_{i_{1},i_{2},\ldots,i_{D}} 
    \varphi_{i_{1}}\pot{1}\kl{x_{1}}  
    \varphi_{i_{2}}\pot{2}\kl{x_{2}} 
    \ldots  
    \varphi_{i_{D}}\pot{D}\kl{x_{D}} ~.
\label{phibasisexpn}   %
\end{eqnarray} 
Re-writing the sum over levels as a sum over basis indices is only this simple if $ m_k(\ell_k) = \ell_k$.  When $ m_k(\ell_k) \ne \ell_k$, the restriction on the basis indices is not the same as the restriction on the levels.   \cite{AvCaV,AvCa}
By equating  \Eq{phibasisexpn} and the Lagrange type function form of the interpolant, 
one obtains an expression for $C_{i_{1},i_{2},\ldots,i_{D}}$ in terms of values of the function 
on the sparse grid. \cite{AvCaV}       In Sec. \ref{sec:trafo},  
we present simpler  ideas for obtaining $ C_{i_{1},i_{2},\ldots,i_{D}}  $.  %
They work only if the sets of points are nested.  

It is common to use a basis of piecewise linear functions, divided into levels, 
that are defined so that the space spanned by the functions in the first $\ell$     
levels is the same as the space spanned by the 
$\ell$-th nodal basis.    \cite{BuGr} %
These functions are called "hierarchical". They are importance-ordered. 
In addition, they have the property that functions in level $\ell$ are equal to zero 
at points in levels $1, 2, \ldots, \ell - 1$.  
We call this the zero-at-points-in-previous-levels (ZAPPL) property.  
When basis functions with the ZAPPL property are used, it is not necessary to sum over levels to determine an interpolant. \cite{halla,BuGr,AvCaV} %
In \cite{VaPf}, the ZAPPL property is called the fundamental property.  
Functions that satisfy the ZAPPL property are also called incremental hierarchical functions. %
 However, in the sparse grid literature they are very often  piecewise linear. \cite{BuGr}
 Ref. \cite{halla} is an important exception to this rule.   Using the ideas of section 2, it is straightforward to make ZAPPL functions from {\textit{any}} global 1-D basis functions.  
 The efficiency  of the approach of Sec. \ref{sec:trafo} relies on using 
basis functions with the ZAPPL property.

\section{ 1-D   ZAPPL  basis   functions}

ZAPPL  functions can be  made  from {\textit{any}}    1-D basis set which is divided    into levels and nested 
sets of points associated with  the levels.  
A general recipe for making  ZAPPL functions  is given in Refs.  \cite{AvCaV,VaPf}.   In those papers the ZAPPL functions are called hierarchical.   %
In this chapter, in 
 level 1, we have  $  \varphi_{1}\pot{k}$, 
in level 2 we have $\kl{ \varphi_{1}\pot{k},  \varphi_{2}\pot{k}   }$, 
in level 3 we have $\kl{    \varphi_{1}\pot{k},  \varphi_{2}\pot{k},  \varphi_{3}\pot{k}   }$, etc.   
Correspondingly, in   level 1, we have  the point  $  r_{1}\pot{k}$, 
in level 2 we have the points $\kl{  r_{1}\pot{k},  r_{2}\pot{k}   }$, 
in level 3 we have $\kl{    r_{1}\pot{k},  r_{2}\pot{k},  r_{3}\pot{k} }$, etc.  
The ZAPPL functions are:
\begin{align}
\tilde{\varphi}_{i_k}\pot{k}\kl{x_k} &= 
\sum_{j_k=1}^{i_k} \tilde{A}\pot{k}_{i_k,j_k} \varphi_{j_k}\pot{k}\kl{x_k},
\label{hier}
\end{align}   
{where   }     $\tilde{A}\pot{k}_{i_k,j_k}   $   is chosen so that $ \tilde{\varphi}_{i_k}\pot{k}\kl{r_{a_k}\pot{k}}=0 \quad  ~~
\forall ~~a_k < i_k \quad \mathrm{and} \quad \tilde{A}\pot{k}_{i_k,i_k}=1$.
Note that 
$
\tilde{\varphi}_{i_{k}}\pot{k}\kl{x_{k}}$ depend on the 
interpolation points and   $
\varphi_{i_{k}}\pot{k}\kl{x_{k}} $ do not.  
The 1-D points must be chosen so that 
$\tilde{\mathbf{B}}\pot{k}$  is  not singular, or near-singular, where    $\tilde{B}\pot{k}_{a_k,i_k} =
\tilde{\varphi}_{i_{k}}\pot{k}\kl{r_{a_k}\pot{k} }$.    
In the rest of this chapter we, for simplicity, omit tildes.   ${B}\pot{k}_{a_k,i_k} $  always means    $\tilde{B}\pot{k}_{a_k,i_k} $   and
${\varphi}_{i_{k}} $  always means  $
\tilde{\varphi}_{i_{k}}$.

The ZAPPL functions defined in \Eq{hier} have the advantageous ZAPPL property, but they may be  smooth and are not the common piecewise linear hierarchical functions. \cite{BuGr}
  To 
interpolate smooth functions it is often better to use smooth basis functions.   
\Eq{hier} can be used to make ZAPPL functions from any importance-ordered basis.  For example, a set of importance-ordered B splines could be used. \cite{VaPf}  
   The prescription of \Eq{hier} can be used regardless of the choice of $ m_k(\ell_k) $, 
   the number of functions in level $ \ell_k $.  %
   In many cases, choosing $ m_k(\ell_k) $ so that it does not increase exponentially with $ \ell_k  $ reduces the cost of
   calculations. 
Of course, it must be possible to choose nested sets of points with  $ m_k(\ell_k)  $ points in level  $ \ell_k  $.  In this chapter, our cost estimates are computed using 
 $ m_k(\ell_k)  =\ell_k $, which means that when $ \ell_k  $ is increased by one, we must add a single new point. One way to do this is to use Leja points. \cite{Leja,NaJa,AvCaX}
\section{Transforming from function values to basis coefficients }
\label{sec:trafo}
In this section, we present our efficient scheme for transforming a vector 
whose elements are values of a function at points on the sparse grid to a vector whose elements  
are the coefficients of a ZAPPL basis expansion of the function.
Everything in this section is valid for {\textit{any}} choice of the 
$ \varphi_{i_{k}}\pot{k}\kl{x_{k}}$ functions and {\textit{any}} choice of the (nested) points.  
We compare our scheme to other transformation methods in the literature. \cite{VaPf,FoTa,Buzz}  %
The method of this section was presented  at the Sparse Grids and Applications conference in Munich.  
After the conference,  David Holzmueller showed that the ideas can be formulated in terms of 
 LU decompositions. \cite{HoPriv}   
Note that if one wishes basis expansion coefficients in a (nested) basis that is not a ZAPPL basis, 
one can use the method of this section to transform from the grid to the ZAPPL coefficients and then 
efficiently (evaluating sums sequentially) transform from the ZAPPL coefficients 
to the coefficients in the desired basis, see Eq. (41) in Ref. \cite{AvCaV}

\subsection{ It appears one needs to invert    $\mathbf{B}$  }
Let $f: \mathbb{R}^D \rightarrow \mathbb{R}$ be a multivariate function.  
Its Smolyak interpolant can be written as in \Eq{phibasisexpn}.
\Eq{phibasisexpn} is similar to the generalised polynomial chaos expansion (GPCE) 
employed when solving stochastic differential equations. \cite{FoTa,XiKa}   
Note, however, that our basis functions are 
not necessarily (weighted) polynomials; they can be chosen to reduce the size of the basis required for the interpolation.  
The goal is to obtain the expansion coefficients $C_\mathbf{i}$ given the values of the function
$f$ at the sparse grid points, i.e,  $f\kl{r_{a_1}\pot{1}, \ldots, r_{a_D}\pot{D}}$. 
Explicit equations for the matrix-vector product required to compute $C_\mathbf{i}$, \Eq{seqsumprun},  
and for the cost, \Eq{nmult2}, are simple if the basis function indices in  \Eq{phibasisexpn} are restricted by $ \norm{\mathbf{i}-\mathbf{1}}_1 \le b   $ and the grid indices are restricted by 
 $ \norm{\mathbf{a}-\mathbf{1}}_1 \le b   $.  
Both these restrictions are inherited from the level restriction in \Eq{bart}.

The most straightforward approach for obtaining $C_\mathbf{i}$ from \Eq{inter} is to solve a system of 
linear  equations
\begin{align}
   \bar{f}\kl{\mathbf{r}_\mathbf{a}} &= 
   f\kl{\mathbf{r}_\mathbf{a}} = \sum_{\norm{\mathbf{i}-\mathbf{1}}_1 \le b} C_\mathbf{i}
   \prod_{k=1}^D \varphi_{i_k}\pot{k}\kl{r_{a_k}\pot{k}} \quad \forall \norm{\mathbf{a}-\mathbf{1}}_1 \le b, \\
   \iff 
   C_\mathbf{i} 
   &= \sum_{\norm{\mathbf{a}-\mathbf{1}}_1 \le b} \klll{\mathbf{B}^{-1}}_{\mathbf{i},\mathbf{a}}
   f\kl{\mathbf{r}_\mathbf{a}} \quad \forall \norm{\mathbf{i}-\mathbf{1}}_1 \le b.
   \label{linsys}
\end{align}
By solving \Eq{linsys}, one obtains all of the coefficients  $C_\mathbf{i}$ from one calculation.  
The elements of the matrix $\mathbf{B}$ are values of the basis functions at the
sparse grid points, i.e., 
\begin{align}
   B_{\substack{a_1, \ldots, a_D \\ i_1, \ldots, i_D}} 
   &= \prod_{k=1}^D \varphi_{i_k}\pot{k}\kl{r_{a_k}\pot{k}} 
   \quad \forall
   \norm{\mathbf{i}-\mathbf{1}}_1 \le b \; \mathrm{and} \;
   \norm{\mathbf{a}-\mathbf{1}}_1 \le b.
   \label{B_eq}
\end{align}
With the chosen restrictions of the indices, the number of sparse grid points (and product basis functions) 
is 
\begin{align}
   N_{\mathrm{sparse}} &= \binom{D + b}{D}.
\end{align}
In  general, to solve the  linear system of equations, $\mathcal{O}\kl{N_{\mathrm{sparse}}^3}$
floating point operations are required.  Directly solving the linear equations in  \Eq{linsys}    or inverting   $\mathbf{B}$
 will therefore, especially for high-dimensional problems, 
require considerable computer time and computer memory. 
Smolyak interpolation thus has the advantage that $N_{\mathrm{sparse}} \ll (b+1)^D$,  
 but the 
disadvantage that it is not simple to determine  $C_\mathbf{i}$.

If both the basis and the point set are tensor products, i.e., the
restrictions imposed on the indices  are $\norm{\mathbf{i}-\mathbf{1}}_\infty \le b$ and
$\norm{\mathbf{a}-\mathbf{1}}_\infty \le b$, then  the entire set of coefficients  $C_\mathbf{i}$ can be easily found 
because one can exploit the fact that 
\begin{align}
   \mathbf{B} &= \mathbf{B}\pot{\mathrm{kron}} = \bigotimes_{k=1}^D \mathbf{B}\pot{k}, 
   \label{kron}
\end{align}
a Kronecker product of small $\kl{b + 1} \times \kl{b + 1}$ matrices
$B_{a_k,i_k}\pot{k}=\varphi_{i_k}\pot{k}\kl{r_{a_k}\pot{k}}$. 
Hence, 
\begin{align}
   \mathbf{B}^{-1} &= \klll{\mathbf{B}\pot{\mathrm{kron}}}^{-1}
   = \bigotimes_{k=1}^D \klll{\mathbf{B}\pot{k}}^{-1}  ~.  
   \label{kroninv}
\end{align}
This has two advantages. First, it is not necessary to invert a large matrix and second, 
it is possible to evaluate the sums in \Eq{linsys} {\textit{sequentially}},
\begin{align}
   C_\mathbf{i} &= 
   \sum_{a_D=1}^{b+1} \kl{\klll{\mathbf{B}\pot{D}}^{-1}}_{i_D,a_D} \ldots 
   \sum_{a_2=1}^{b+1} \kl{\klll{\mathbf{B}\pot{2}}^{-1}}_{i_2,a_2}
   \sum_{a_1=1}^{b+1} \kl{\klll{\mathbf{B}\pot{1}}^{-1}}_{i_1,a_1}
   f\kl{\mathbf{r}_\mathbf{a}},
   \label{seqsum}
\end{align}
as is frequently done in chemical physics to  transform  between a grid and a basis.  \cite{LiCa} %
Not exploiting the Kronecker product structure of \Eq{kron} and solving the
linear system directly would require 
$\mathcal{O}\kl{\kl{b + 1}^{3D}}=\mathcal{O}\kl{N_\mathrm{full}^3}$
operations, whereas the sequential summation approach of \Eq{seqsum} requires
only 
\begin{align}
   \mathcal{O}\kl{D\kl{b + 1}^{D + 1}}
   &= \mathcal{O}\kl{D\kl{b+1}N_\mathrm{full}} 
   \label{fullscal}
\end{align}
operations, where $N_\mathrm{full}=\kl{b + 1}^D$. 
	Each sum in Eq. (13) can be thought of as calculating $(b+1)^{D-1} $ MVPs for a matrix of size $ (b+1)$.  
      As each of these  MVPs requires     $(b+1)^{2}$ multiplications and there are     $(b+1)^{D-1}$ of them, the total cost scales as $(b+1)^{D+1}$.   
 \Eq{fullscal} is the cost of computing the entire set of $C_\mathbf{i}$ coefficients.

Although $\mathbf{B}$  in \Eq{linsys} is {\textit{not}} a tensor product, we show in the next subsection that 
the sparse grid to basis transformation can nevertheless be done sequentially, when ZAPPL basis functions are employed.
The numerical cost then scales as
\begin{align}
   &\mathcal{O}\kl{D \klll{\frac{b}{D+1} + 1} N_\mathrm{sparse}}~.
   \label{sparsescal}
\end{align}
Again, this is the cost of computing the entire set of $C_\mathbf{i}$ coefficients.  
Because in both cases sums are evaluated sequentially,  \Eq{sparsescal} 
and the right side of \Eq{fullscal} have the same structure    and there is a factor of $D$ in both,  but in  
\Eq{sparsescal}, $N_\mathrm{full}$ is replaced with $N_\mathrm{sparse}$ 
and $\kl{b + 1}$ is replaced with $\frac{b}{D+1}+1$.  %

\subsection{ZAPPL functions and sequential summation obviate the need to invert  $\mathbf{B}$}
\label{sec:seqsum}
We begin by sorting  the  products of ZAPPL functions that form the tensor product basis  into two groups:
the $  N_\mathrm{sparse}$   functions   
 that are included in the sparse basis go into a group  labelled "retained"; and  the  $ {N_\mathrm{full} - N_\mathrm{sparse}}$  functions   
 that are excluded from the sparse basis go into a group labelled  "discarded".
Let $\mathbf{C} \in \mathbb{R}^{N_\mathrm{full} \times N_\mathrm{sparse}}$ be a 
chopping matrix which  is  an $N_\mathrm{full} \times N_\mathrm{full}$ 
identity matrix from which the  columns  for  the excluded functions  have been
deleted. %
In terms of   $\mathbf{C} $,   $\mathbf{B}$    %
can be written  as
\begin{align}
   \mathbf{B} &= \mathbf{C}^T \mathbf{B}\pot{\mathrm{kron}} \mathbf{C} \\
   \iff
   \mathbf{B}^{-1} &= \klll{\mathbf{C}^T \mathbf{B}\pot{\mathrm{kron}} \mathbf{C}}^{-1}.
   \label{binv}
\end{align}
It is far from obvious that one can calculate elements of the matrix on the right 
side of \Eq{binv} by inverting small matrices for each coordinate and  then  do  the  sums    in an equation like   %
\Eq{linsys} sequentially so that one obtains an equation similar to \Eq{seqsum}.    %
Both are necessary, if one is to find an inexpensive method for computing   $   C_\mathbf{i} $  on the LHS of \Eq{linsys}.   %
However, if excluded functions are zero at retained points then
\begin{align}
	\klll{\mathbf{C}^T \mathbf{B}\pot{\mathrm{kron}} \mathbf{C}}^{-1}
	&=
	\mathbf{C}^T \klll{\mathbf{B}\pot{\mathrm{kron}}}^{-1} \mathbf{C} ~.
	\label{trick}
\end{align}
In  words, we can interchange the operations of chopping and inverting. 
The  matrix in the  matrix-vector product in \Eq{linsys} is the inverse of the  retained block, but because the inverse of the retained block can be replaced by 
a block of the inverse of ${\mathbf{B}\pot{\mathrm{kron}}}  $ (due to \Eq{trick}) it is possible do 
sums sequentially; see \Eq{seqsumprun}.
 \Eq{trick} will be  satisfied for any order of the tensor product  basis functions if it is satisfied for one order. By  
 block Gaussian elimination,  \cite{block}   %
 it is  simple  to prove  that  \Eq{trick} is
 correct if the tensor product  basis functions are ordered so that
 the 
  top left block of the reordered  $ \mathbf{B}\pot{\mathrm{kron}} $   is $ 	\klll{\mathbf{C}^T \mathbf{B}\pot{\mathrm{kron}} \mathbf{C}} $  
  and   %
the top right block of the reordered    $   \mathbf{B}\pot{\mathrm{kron}}   $  is zero.  %
 We have done calculations with the  %
 $
 \lVert \mathbf{i} -\mathbf{1} \rVert \le b$ pruning condition and in this case it is easy to show that if the 1-D functions are ZAPPL functions then functions excluded from the tensor product basis are zero at retained points and thus \Eq{trick} is satisfied.

We are now able to evaluate the sum in \Eq{linsys} sequentially,  
and when using a simple pruning condition obtain explicit equations for the upper limits of the sums.  %
For example, for the $ \norm{\mathbf{a}-\mathbf{1}}_1 \le b$ pruning condition, one has,
\begin{align}
   C_\mathbf{i} &= \sum_{\norm{\mathbf{a}-\mathbf{1}}_1 \le b} 
   \kl{\mathbf{C}^T \klll{\mathbf{B}\pot{\mathrm{kron}}}^{-1} \mathbf{C}}_{\mathbf{i},\mathbf{a}} 
   f\kl{\mathbf{r}_\mathbf{a}} 
   \nonumber \\
   &= \sum_{\norm{\mathbf{a}-\mathbf{1}}_1 \le b} 
   \kl{\klll{\mathbf{B}\pot{\mathrm{kron}}}^{-1}}_{\mathbf{i},\mathbf{a}}
   f\kl{\mathbf{r}_\mathbf{a}} \quad \forall \norm{\mathbf{i}-\mathbf{1}}_1 \le b \nonumber \\ 
   C_\mathbf{i} &= 
   \sum_{a_D=1}^{a_D\pot{\mathrm{max}}} \kl{\klll{\mathbf{B}\pot{D}}^{-1}}_{i_D,a_D} \ldots 
   \sum_{a_2=1}^{a_2\pot{\mathrm{max}}} \kl{\klll{\mathbf{B}\pot{2}}^{-1}}_{i_2,a_2}
   \sum_{a_1=1}^{a_1\pot{\mathrm{max}}} \kl{\klll{\mathbf{B}\pot{1}}^{-1}}_{i_1,a_1}
   f\kl{\mathbf{r}_\mathbf{a}} ~.   
   \label{seqsumprun}
\end{align}
If the $ [{\mathbf{B}\pot{k}} ]^{-1}$ matrices were not lower triangular, 
then the range of possible $i_k$ values would be limited by %
\begin{align}
   i_k\pot{\mathrm{max}} 
   &=
   b - \sum_{k' = 1}^{k-1}  \kl{i_{k'} - 1} + 1 ~,
\end{align}
and the upper limits on the sums over $ a_k $ would be   
\begin{align}
a_k\pot{\mathrm{max}} 
   &=
   b - \sum_{k'= k + 1}^{D} \kl{a_{k'} - 1} + 1~.
\end{align}
However, because ${\mathbf{B}\pot{k}}$ is lower triangular, 
$[{\mathbf{B}\pot{k}} ]^{-1}$ is also lower triangular and therefore the upper limits of the sums are
\begin{align}
   a_k\pot{\mathrm{max}} 
   &= i_k.
\end{align}
and 
\begin{align}
   i_k\pot{\mathrm{max}} 
   &= 
	b - \sum_{k' = 1}^{k-1}  \kl{i_{k'} - 1} 
	- \sum_{k'= k + 1}^{D} \kl{a_{k'} - 1} + 1.  
\end{align}
Related  sequential summation techniques were previously used:  
with a basis restricted by ${\norm{\mathbf{j}-\mathbf{1}}_1 \le b} $ in Ref. \cite{WaCa}; 
for Smolyak quadrature in Refs.  \cite{AvCa,AvCaIII};  %
for collocation in Ref.   \cite{AvCaVIII};
and for collocation with a hierarchical (ZAPPL) basis in Ref. \cite{AvCaV,AvCaIX,AvCaVI,emil,WoCaIV}.
Note that in Refs. \cite{ AvCa,AvCaV,AvCaIX,AvCaVI,emil} a similar  sequential summation idea is used with sparse grids made with  either $ m_k(\ell_k) \ne \ell$ or a pruning condition more general than $
{\norm{\mathbf{i}-\mathbf{1}}_1 \le b}  $.
The unidirectional principle also exploits the sequential evaluation of sums.  
\cite{balder,feuersaenger,dirk,zeiser}

The numerical cost  associated with the prescription in \Eq{seqsumprun} can be
estimated  by counting the number of multiplications required 
for each of the $D$ sequential summations. 
The number of multiplications required for the sum   over    
$a_1$ is  
\begin{align}    
   \frac{1}{D}N_\mathrm{mult}
   &=
   \sum_{a_D     = 1}    ^{b + 1} \;
   \sum_{a_{D-1} = 1}    ^{b - \kl{a_D - 1} + 1} \ldots 
   \sum_{a_2     = 1}    ^{b - \sum_{k'=3}^D \kl{a_{k'} - 1} + 1} \;
   \sum_{i_1     = 1}    ^{b - \sum_{k'=2}^D \kl{a_{k'} - 1} + 1} \;
   \sum_{a_1     = 1}    ^{i_1} 1 \nonumber \\
   &=
   \sum_{a_D     = 0}    ^{b}
   \sum_{a_{D-1} = 0}    ^{a_D} \ldots 
   \sum_{a_2     = 0}    ^{a_3} 
   \sum_{i_1     = 0}    ^{a_2} 
   \sum_{a_1     = 0}    ^{i_1} 1 \nonumber \\
   &=
   \sum_{a_D     = 0}    ^{b}
   \sum_{a_{D-1} = 0}    ^{a_D} \ldots 
   \sum_{a_2     = 0}    ^{a_3} 
   \sum_{i_1     = 0}    ^{a_2} \kl{i_1 + 1}.
   \label{nmult}
\end{align}
The expression in \Eq{nmult} can be evaluated analytically \cite{BuKa} yielding 
\begin{align}
   N_\mathrm{mult} &= D\kl{\frac{b}{D+1} + 1}N_\mathrm{sparse},
   \label{nmult2}
\end{align}
which is  the scaling given in \Eq{sparsescal}.   \Eq{sparsescal} is the cost of computing all the   $C_\mathbf{i}$ coefficients.  
The transformation of \Eq{seqsumprun} was used in Ref. \cite{WoCaIV,emil}.
Its cost depends directly on the number of points on the sparse grid.

\subsection{Comparison with other methods}

In the field of uncertainty quantification, \cite{uq1,uq2}  
 it is often desirable to transform 
from a vector of values of a function on a  sparse grid to a vector of coefficients of a  so-called generalised polynomial chaos expansion (GPCE), for example,  to facilitate
obtaining stochastic quantities of interest.
The GPCE is a special case of the sparse grid interpolation  of
 \Eq{bart}, where in the GPCE case  the 
$\varphi_{i_k}\pot{k}\kl{x_k}$ are (weighted) polynomials.
Several different methods are used to obtain the coefficients.
Some methods for computing basis expansion coefficients build a $\mathbf{B}$ matrix 
whose rows and columns each have a level label and a label identifying a basis function within that level. \cite{VaPf}    %
Although  $\varphi_{i_k}\pot{k}\kl{x_k}$ can be chosen to make it possible to exploit 
the sparsity of the corresponding  $\mathbf{B}$, the sums required to compute the coefficients are costly.  %
The key advantage of our algorithm is the sequential evaluation of the sums.  
Algorithm 1 of Ref. \cite{VaPf} scales as the square of the number of functions in the basis.  
The scaling of  \Eq{nmult2} is much better.  
Other methods for computing basis expansion coefficients manipulate separately 
the grids that together compose the sparse grid. An algorithm of this kind is 
used in Refs. \cite{FoTa,Buzz}. They solve linear systems for all the 
tensor product grids associated with the levels $\bm{\ell}$ that satisfy
$\norm{\bm{\ell}-\mathbf{1}}_1 \le b$.
For one of the tensor product grids, they  write 
\begin{align}
   f\kl{\mathbf{r}_\mathbf{\tilde{a}}} &=                 
   \sum_{p_k \le \ell_k \; ~ \forall k } \;
   \beta^{\bm{\ell}}_{\mathbf{p}}    %
 ~ L^{\bm{\ell}}_{\mathbf{p}}\kl{\mathbf{r}_\mathbf{\tilde{a}}}     ~,   
   \label{tamel}            %
\end{align}
where $ L^{\bm{\ell}}_{\mathbf{p}}   $ is a basis function for the tensor product grid 
labelled by $\bm{\ell}$,  
and solve for $ \beta^{\bm{\ell}}_{\mathbf{p}}   $  for each tensor product grid separately.   
Coefficients with the same $\mathbf{p}$ for different grids must 
be combined to determine $C_\mathbf{i}$.   
$   \mathbf{r}_{\mathbf{\tilde{a}}  }$ is a point on the tensor product grid associated with $\bm{\ell}$;  
the union of the grids whose points are denoted 
 $   \mathbf{r}_{\mathbf{\tilde{a}}  }$ 
 is the sparse grid on which the points are denoted by $   \mathbf{r}_{\mathbf{{a}}  }$.

\begin{figure}[H]
   \caption{Scaling of floating point operations as a function of $D$, $b=4$ (i.e. 5 points per coordinate)   }    
   \label{fig:fig1}
   \includegraphics[angle=0,width=1.0\linewidth]{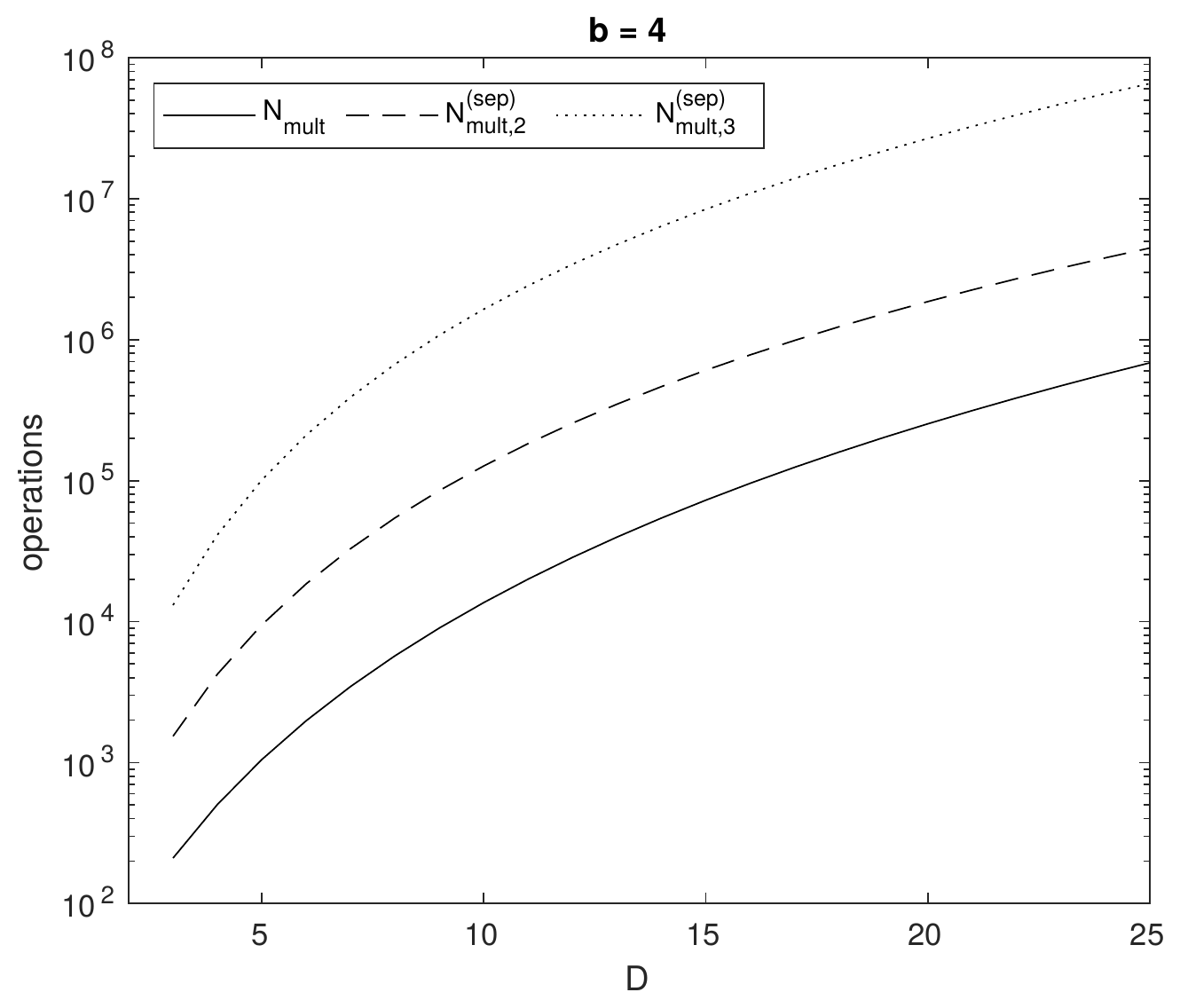}
\end{figure}

\begin{figure}[H]
   \caption{Scaling of floating point operations as a function of $D$, $b=9$ (i.e. 10 points per coordinate) }  
   \label{fig:fig2}
   \includegraphics[angle=0,width=1.0\linewidth]{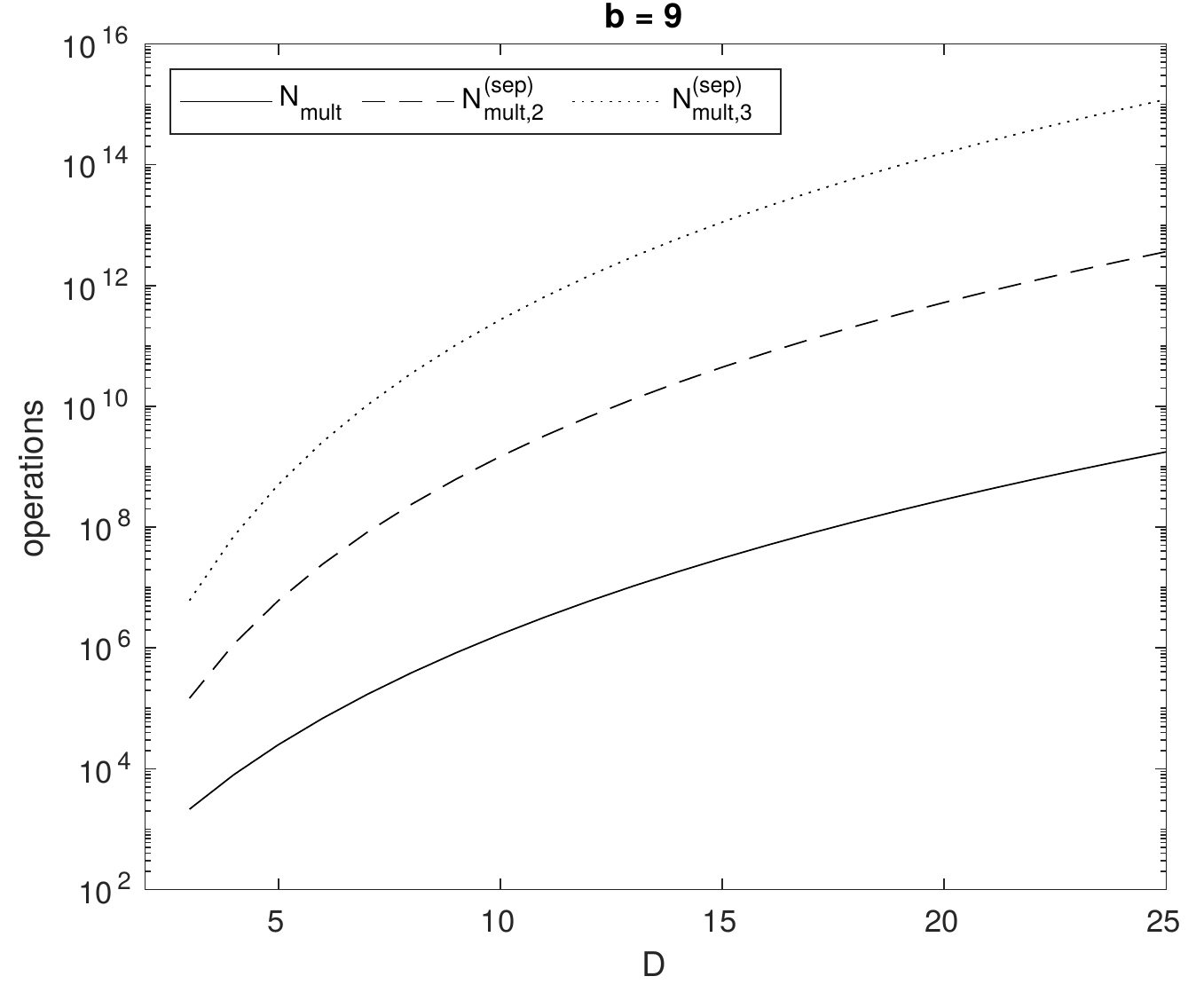}
\end{figure}

\begin{figure}[H]
	\caption{Scaling of floating point operations as a function of $D$, $b=14$ (i.e. 15 points per coordinate) } 
	\label{fig:fig3}
	\includegraphics[angle=0,width=1.0\linewidth]{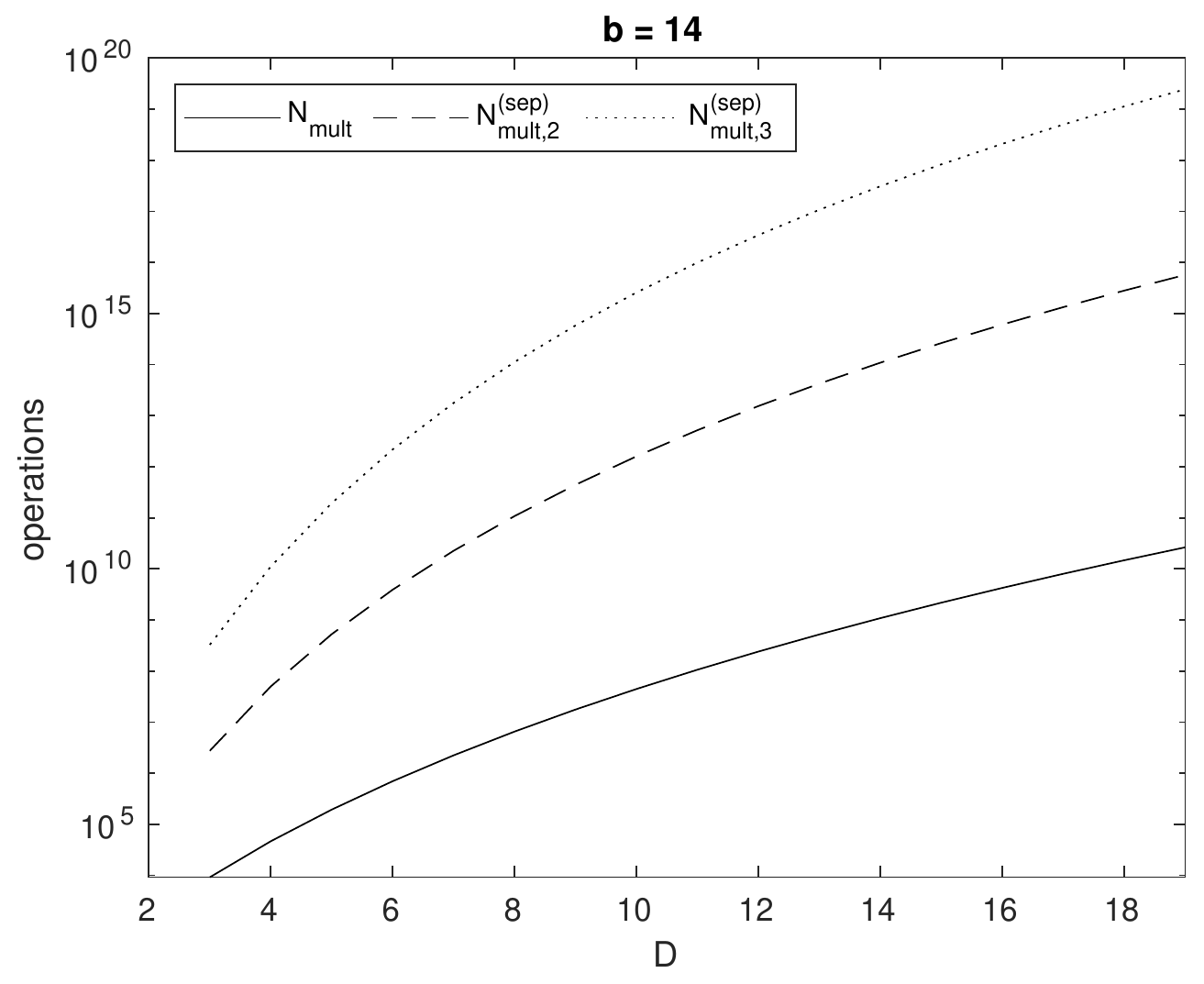}  %
\end{figure}

The total cost  of this separate grids method    is the sum of the costs of matrix-vector products for the 
 grids, $N_{\mathrm{mult},2}\pot{\mathrm{sep}}$,
and the cost of inverting matrices for the  grids,
$N_{\mathrm{mult},3}\pot{\mathrm{sep}}$, 
 \begin{align}
    \mathrm{cost} &=   
   N_{\mathrm{mult},2}\pot{\mathrm{sep}}   +   N_{\mathrm{mult},3}\pot{\mathrm{sep}}~,
  \label{cost}
 \end{align}
 where 
 \begin{align}
 N_{\mathrm{mult},n}\pot{\mathrm{sep}}  %
 &=
 \sum_{\norm{\bm{\ell}-\mathbf{1}}_1 \le b} \;
 \prod_{k=1}^D \ell_k^n.
 \label{nmultgpce}
 \end{align}
 In this equation, we assume  that the number of points for coordinate $x_k$ in level 
 $\ell_k$ is equal to $\ell_k$ and that
 the cost of the matrix-vector product 
 required to obtain $ \beta^{\bm{\ell}}_{\mathbf{p}} $ is $ ( \prod_{k=1}^D \ell_k)^2  $.
If the matrix whose elements are 
$  L^{\bm{\ell}}_{\mathbf{p}}\kl{\mathbf{r}_\mathbf{\tilde{a}}}  $
is orthogonal,  
 and there is therefore no need to invert it, 
$ N_{\mathrm{mult},2}\pot{\mathrm{sep}}   $ is the total cost.   
On the other hand, when  $  L^{\bm{\ell}}_{\mathbf{p}}\kl{\mathbf{r}_\mathbf{\tilde{a}}}  $  is not  orthogonal, 
 then $  N_{\mathrm{mult},3}\pot{\mathrm{sep}} $ must be included in \Eq{cost}.
It is also possible to exploit the tensor product character of the grids.  
After some algebraic transformations, \Eq{nmultgpce} can be cast
into  a form similar to \Eq{nmult}.
However, it does not seem to be possible to find an equation in closed form for $N_{\mathrm{mult},n}\pot{\mathrm{sep}}$
for a general $D$. \cite{BuKa}   %
We  thus evaluate \Eq{nmultgpce} numerically, to compare the cost of this separate-grids method and the approach in 
Sect. \ref{sec:seqsum}.  
We plot the numbers of operations, in \Eq{nmult2} and \Eq{nmultgpce},
as a function of $D$ for the threshold parameters $b=4$, $b=9$, and $ b=14 $  
in Figs. \ref{fig:fig1}, \ref{fig:fig2}, and \ref{fig:fig3}, respectively.   %
For  Figure 1, Figure 2, and Figure 3,  the columns that are deleted  are those for which  $\norm{\mathbf{i}-\mathbf{1}}_1 > b$.

When $b=9$ and $ b=14 $ the sequential summation method of Sect. \ref{sec:seqsum} is orders of magnitude 
less costly than the separable grids method.

\section{Conclusion}

In this chapter, we explain  that by doing sums sequentially it is possible to efficiently 
obtain  expansion coefficients in a basis of products of 1-D global functions  %
from values of a function on a sparse grid. 
Such coefficients are needed in many contexts. They are needed whenever one wishes to make
an interpolant from a sparse grid and global 1-D functions.     
For example, in Sect. 2.6 of Ref. \cite{dakota} the coefficients, called polynomial chaos expansion coefficients,
are determined by evaluating integrals with  Smolyak quadratures.    
In Sect. 2.7 of Ref. \cite{dakota}, the coefficients are obtained by solving a system of linear equations 
whose size is equal to the number of sparse grid points.    
If both the basis and the grid are tensor products, it is well established that basis coefficients 
can be computed from function values by doing sums sequentially. See for example, Table 10.1 of Ref. \cite{boyd}.  
By doing sums sequentially, one obviates the need to loop simultaneously over all the indices.
It is not well established that a similar sequential summation idea can be used when the grid 
and the basis are built from a general sparse grid recipe and nested sets of 
1-D grids and global 1-D bases.     %
In Ref.  \cite{AvCaV}, the original sparse grid interpolation, in terms of 
products of differences of 1-D interpolation operators, was formulated using 
sequential sums.   
In this chapter, we have shown that it is also possible to interpolate, i.e., 
to compute basis expansion coefficients, using sequential sums, without writing the 
interpolant in terms of products of differences of 1-D interpolants. 
The sequential summation approach is considerably less costly than established methods.

\section{Relation with the Chapter of David Holzmueller and Dirk Pflueger}

We have  read a preliminary version of the  chapter by
Holzmueller and  Pflueger and wish to clarify the relationship between their contribution and ours. 
  First, 
  in many of our equations 
   we use a special and simple restriction to determine which functions (points) are included in the sparse basis (grid): $
   \lVert \mathbf{i} - \mathbf{1} \rVert \le b$. 
    This enables to write out the sums in \Eq{seqsumprun} explicitly and to impose the restriction by using the appropriate upper limits on the sums.  
Holzmueller  and  Pflueger do not explain in detail how they evaluate matrix vector products, but appear to use index lists and 
their formalism allows  one to use much more general restrictions.  
  Ideas similar to those in this paper  
can be used with the restriction $g_1(i_1-1) + g_2(i_2-1) + \cdots  + g_D(i_D-1) \le b$, 
where $ g_c(i_c-1), c= 1, \cdots D $, 
 is a monotonically increasing function.\cite{AvCaV,emil}. 
  In this case it is also possible to derive equations for the upper
  limits on the sums in \Eq{seqsumprun}.  Much more general restrictions can be used if one is willing to forgo deriving equations 
  for the upper limits and instead uses index lists.\cite{ourfirst,oursecond}      
Second, in this paper we have one function (point) per level.  Holzmueller  and  Pflueger have no such constraint.
  However, this restriction can be lifted.\cite{AvCaV}.  
  Third,  we have not given a recipe for choosing 1-D basis functions and dividing the tensor product functions obtained from  them  into two groups, labelled "retained" and "discarded",  so that all the discarded functions are zero at the retained points.  This is simple if the 1-D functions are ZAPPL functions and $ all $
  the functions that satisfy  $g_1(i_1-1) + g_2(i_2-1) + \cdots  + g_D(i_D-1) \le b$ are retained ($ g_c(i_c-1) $ is a monotonically increasing function).   It is  in general not simple if $ some $ of the 
 functions that satisfy  $g_1(i_1-1) + g_2(i_2-1) + \cdots  + g_D(i_D-1) \le b$ are excluded, i.e., if there are holes.   
 An easy, and often inexpensive, way to avoid holes is to plug them by adding functions to the retained basis. 
       When using spatially localised basis functions, it is sometimes advantageous   to discard functions 
       in regions in which the function being interpolated is smooth, i.e. to introduce holes.    In the sparse grid literature this is 
       known as spatial adaptivity.   When using global basis functions,  often best for representing a smooth function, spatial adaptivity is less useful.

\section*{Acknowledgements}
Research reported in this article was funded by  The Natural Sciences and Engineering Research Council of Canada.  
We thank David Holzmueller for sending us reference \cite{Buzz}, for emails about his LU perspective.  We are grateful to both David Holzmueller and Dirk Pflueger for discussions.

\pagebreak

\end{document}